%% file: twisted-main.tex
\documentclass[reqno,11pt]{amsart}

\usepackage{mathrsfs}
\usepackage{graphicx}
\usepackage{amscd}
\usepackage{amsmath}
\usepackage{amsthm}
\usepackage{amsfonts}
\usepackage{amssymb}
\usepackage{amsxtra}
\usepackage{xspace}
\usepackage{array}
\usepackage{cite}
\usepackage{color}

\theoremstyle{plain}
\newtheorem{theorem}{Theorem}
\newtheorem{lemma}{Lemma}
\newtheorem{corollary}{Corollary}
\newtheorem{proposition}{Proposition}

\theoremstyle{definition}

\newtheorem{remark}{Remark}

\numberwithin{equation}{section}

\newcommand{\be}{\begin{enumerate}}
\newcommand{\ee}{\end{enumerate}}
\newcommand{\beq}{\begin{equation}}
\newcommand{\eeq}{\end{equation}}
\newcommand{\bprop}{\begin{proposition}}
\newcommand{\eprop}{\end{proposition}}
\newcommand{\complex}{\mathbb{C}}
\newcommand{\reals}{\mathbb{R}}

\newcommand{\integers}{\mathbb{Z}}
\newcommand{\pfbegin}{\noindent {\em Proof:} }

\DeclareMathOperator{\htt}{ht} 
\DeclareMathOperator{\ct}{ct} 
\DeclareMathOperator{\mult}{mult}

\newcommand{\s}{\scriptscriptstyle}

\newcommand{\form}[2][]{\langle #1,#2 \rangle}
\newcommand{\poch}[2][]{(#1;\, #2)_{\s \infty}}

\newcommand{\alolost}[1][s,t,q]{\ensuremath{a^{\s \Lambda_0}_{\s \Lambda_0}(#1)}\xspace}
\newcommand{\alolo}[1][t,q]{\ensuremath{a^{\s \Lambda_0}_{\s \Lambda_0}(#1)}\xspace}
\newcommand{\alm}[1][t,q]{a^{\s \lambda}_{\s \mu}(#1)} 
 
\newcommand{\tko}[1][]{\mathcal{P}(#1; \,t)}
\newcommand{\stko}[1][]{\mathcal{P}(#1; \,s,t)}
\newcommand{\kost}[1][\mu]{K_{\lambda #1}(t)}
\newcommand{\kostkast}[1][]{K_{#1}(s,t)}

\newcommand{\sch}[1][\lambda]{\chi_{{}_{#1}}}

\newcommand{\kma}{\mathfrak{g}}
\newcommand{\kl}{\mathfrak{m}}
\newcommand{\nilpf}{\mathfrak{n}}

\newcommand{\ph}{\mathfrak{s}}
\newcommand{\csa}{\mathfrak{h}}
\newcommand{\bas}{L(\Lambda_0)}

\newcommand{\chkeim}{\chke^{im}}
\newcommand{\cato}{\mathcal{O}}

\newcommand{\kmaf}{\overset{\circ}{\mathfrak{g}}}

\newcommand{\chke}{\hat{\boldsymbol{\mu}}}
\newcommand{\spchke}{\boldsymbol{\varDelta}}

\newcommand{\Qf}{\overset{\circ}{Q}}

\newcommand{\tsf}{$t$-string function\xspace}
\newcommand{\tsfs}{$t$-string functions\xspace}

\newcommand{\ex}{\mathbb{E}}

\newcommand{\rof}{\overset{\circ}{\Delta}}

\newcommand{\ro}{\Delta}
\newcommand{\rp}{\Delta_+}

\newcommand{\rre}{\Delta^{re}}
\newcommand{\rrep}{\Delta^{re}_+}
\newcommand{\rim}{\Delta^{im}}
\newcommand{\rimp}{\Delta^{im}_{+}}

\newcommand{\ccheckc}{(C_l^\vee, C_l)}
\newcommand{\cc}{\Phi}

\newcommand{\dprod}{\displaystyle\prod}

\begin{document}

\title[]{The $t$-analog of the basic string function for twisted affine Kac-Moody algebras}
\author{Sachin S. Sharma and Sankaran Viswanath}
\address{The Institute of Mathematical Sciences\\
CIT campus, Taramani\\
Chennai 600113, India}
\email{sachin@imsc.res.in, svis@imsc.res.in}

\begin{abstract}
We study Lusztig's $t$-analog of weight multiplicities associated to level one representations of twisted affine Kac-Moody algebras.  An explicit closed form expression is obtained for the corresponding $t$-string function using constant term identities of Macdonald and Cherednik. The closed form involves the generalized exponents of the graded pieces of the twisted affine algebra, considered as modules for the underlying finite dimensional simple Lie algebra. This extends previous work on 
level 1 $t$-string functions for the untwisted simply-laced affine Kac-Moody algebras.

\end{abstract}
\maketitle

\input{twisted-tsf.tex}

\input{twisted-main.bbl}
\end{document}

%% file: twisted-tsf.tex
Let $\kma$ be an irreducible affine Kac-Moody algebra (over $\complex$). If $\kma$ is untwisted, it can be realized as a central extension of a loop algebra over a finite dimensional simple Lie algebra. Twisted affines are realized as fixed point subalgebras of untwisted affine Lie algebras under finite groups of automorphisms. The category $\cato^{\mathrm{int}}(\kma)$ of integrable $\kma$-modules in category $\cato$ is semisimple, and its simple objects are highest weight representations $L(\lambda)$ indexed by dominant integral weights $\lambda$. The formal character of $L(\lambda)$ is given explicitly by the Weyl-Kac character formula. To understand the structure of the module $L(\lambda)$, one studies the generating function 
$$\alm[q] :=  \sum_{k \geq 0} \dim (L(\lambda)_{\mu - k\delta})\; q^k$$
of weight multiplicities along $\delta$-strings through dominant weights $\mu$ of $L(\lambda)$. Here $\delta$ is the null root of $\kma$. One usually also assumes that $\mu$ is {\em maximal}, i.e., that $\mu + \delta$ is not a weight of $L(\lambda)$. We will call the $\alm[q]$ {\em string functions} of the module $L(\lambda)$ (mildly departing from standard convention \cite{kac}).

Among the irreducible modules in  $\cato^{\mathrm{int}}(\kma)$, the so called {\em basic representation} $L(\Lambda_0)$ can be singled out for the unique and important role it plays in the theory. It is the simplest non-trivial representation in  $\cato^{\mathrm{int}}(\kma)$, and has many explicit realizations in terms of vertex operators \cite{kac-kazhdan-lepow-wilson, frenkel-kac, kac-peterson-112constr}. If  $\kma$ is an untwisted simply-laced affine Lie algebra, or if $\kma$ is twisted, then the basic representation is an irreducible $\kma$-module of level one, and all level one simple modules can be obtained from it by the action of Dynkin diagram automorphisms of $\kma$, and by tensoring with one-dimensional $\kma$-modules.

Further, for such $\kma$,  the basis representation admits a unique  string function $\alolo[q]$, which we will refer to as the {\em basic string function}. This has a well known \cite{kacpeterson, kac}, closed form expression as an infinite product (equation \eqref{bsf}). For the remaining affine Kac-Moody algebras, i.e, the untwisted ones of types $B, C, F, G$, there are multiple inequivalent string functions of level 1. 

Now, let $t$ be an indeterminate.  The $t$-analog $\tko[\beta]$ of the Kostant partition function is defined to be the coefficient of $e^\beta$ in the product 
$$\prod_{\alpha \in \rp(\kma)} (1- te^\alpha)^{- \mult \alpha}$$
where $\rp(\kma)$ is the set of positive roots of $\kma$ and $\mult \alpha := \dim \kma_\alpha$. In this article, we will be concerned with Lusztig's $t$-analog of weight multiplicity, or (affine) Kostka-Foulkes polynomial $\kost$. Given a dominant integral weight $\lambda$, and a weight $\mu$ of $L(\lambda)$, $\kost$ is defined to be : 
$$\kost := \sum_{w \in W} \epsilon(w) \, \tko[w(\lambda + \rho) - (\mu + \rho)]$$ 
where $\epsilon$ is the sign character of the Weyl group $W$ of $\kma$. These polynomials have several important properties: (a) they are the transition coefficients between the bases of affine Hall-Littlewood functions and the formal characters of simple modules \cite{svis-kfp} (b) they have non-negative integral coefficients \cite{slofstra} and (c) $K_{\lambda \mu}(1) = \dim L(\lambda)_\mu$. In view of (c), it is natural to consider the $t$-analog of the string function, defined by $$\alm := \sum_{k \geq 0} K_{\lambda, \mu-k\delta}(t) \; q^k $$ We call these the $t$-string functions. It was shown in \cite{svis-kfp} that the $\alm$ are closely related to the constant term identities arising in the theory of Macdonald polynomials. For instance, when $\kma$ is an untwisted simply-laced affine Kac-Moody algebra, Cherednik's Macdonald-Mehta constant term identity \cite{dmmc} allows us to compute the $t$-string function $\alolo[t,q]$ of the basic representation, in closed form\footnote{We point out an erratum in \cite{svis-kfp}. Theorems 2, 3 and Corollary 3 there are missing the hypothesis that $\kma$ must be simply laced.}. More precisely \cite[corollary 2]{svis-kfp}:

\begin{theorem}\label{ade}
Let $\kma$ be one of the simply laced untwisted affine Lie algebras $A_l^{(1)}, D_l^{(1)}, E_l^{(1)}$. Then
$\alolo = \dprod_{i=1}^l \dprod_{n=1}^\infty (1-t^{e_i+1} q^n)^{-1}$
where $e_i \,(1 \leq i \leq l)$ are the exponents of the underlying finite
    dimensional simple Lie algebra ($=A_l, D_l$ or $E_l$ respectively).
\end{theorem}

The goal of the present work is to prove the twisted version of this theorem.  Our main result  gives 
a closed form expression for the basic $t$-string function of the twisted affine Kac-Moody algebras. 

We recall that every  affine Kac-Moody algebra $\kma$ admits a natural $\integers$-grading into finite dimensional subspaces $\kma = \oplus_{j \in \integers} \kma_j$, where $\kma_0$ is a reductive Lie algebra, and the $\kma_j$ are irreducible $\kma_0$-modules. For $\kma$ twisted, the basic $t$-string function (theorem \ref{mainthm}) now involves the generalized exponents of the $\kma_j$ ($j>0$), viewed as representations of the semisimple part of $\kma_0$. When $\kma$ is untwisted simply-laced, all the $\kma_j$ are isomorphic to the adjoint representation of the underlying finite dimensional simple Lie algebra. Since in this case, the generalized exponents are just the usual exponents, our theorem is in fact a generalization of theorem \ref{ade}.

The paper is organized as follows. The required preliminaries and the statement of the main theorem appear in 
\S\ref{first}. The proof for the case $\kma \neq A_{2l}^{(2)}$ is in \S \ref{pfa}. When $\kma = A_{2l}^{(2)}$, as an auxiliary step, we first consider the non-reduced affine root system $\ccheckc$. We show that on passing to an appropriate limit, the Macdonald constant term identity for this root system can be used to derive the Macdonald-Mehta type constant term identity for $A_{2l}^{(2)}$. This is used to prove the theorem in this case. As an offshoot, we define a {\em two-variable} generalization of the Kostka-Foulkes polynomials for $A_{2l}^{(2)}$, and show that the two variable polynomials associated to the basic representation have non-negative integral coeffients. These details appear in \S\ref{pfb}. The main results of this article were announced in \cite{announcement}.

\section{The main theorem}\label{first}
\subsection{Preliminaries}
For use in stating our main theorem, we recall the definition of generalized exponents.
 Let $\kma$ be a finite dimensional simple Lie algebra and $V =V(\lambda)$ be the irreducible finite dimensional 
$\kma$-module with highest weight $\lambda$. Fix a triangular decomposition $\kma = \nilpf^- \oplus \csa \oplus \nilpf^+$, 
and let $E \in \nilpf^+$ be a regular (principal) nilpotent
element.  Let $V_0$ denote the zero weight space of $V$, and define the Brylinski-Kostant filtration \cite{rkg,kostant3d} of $V_0$ via
 $\mathcal{F}^{(p)}(V_0) := \ker\,(E^{p}) \cap V_0$ for $p \geq 0$. Then the multiset $\ex(V)$ of generalized exponents of $V$ is defined via the relation: 
$$ \sum_{p \geq 0} \dim(\mathcal{F}^{(p+1)}(V_0)/\mathcal{F}^{(p)}(V_0)) \, t^p = \sum_{k \in \ex(V)} t^k.$$ 

Next, given a Kac-Moody algebra $\kma$ of finite or affine type, we let $\ro(\kma)$, $\rp(\kma)$, $\rre(\kma)$, $\rim(\kma)$ denote the sets of roots, positive roots, real roots and imaginary roots respectively, and let $\rrep(\kma) := \rre(\kma) \cap \rp(\kma)$, $\rimp(\kma) := \rim(\kma) \cap \rp(\kma)$.
The {\em Cherednik kernel} $\chke$ of $\kma$ is the product 
\beq \label{which}
\chke :=\prod_{\alpha \in \rrep(\kma)} \frac{1-e^{-\alpha}}{1-te^{-\alpha}}
\eeq
When $\kma$ is affine, we denote the corresponding product over the imaginary positive roots by 
$\chkeim := \prod_{\alpha \in \rimp(\kma)} \left( \frac{1-e^{-\alpha}}{1-te^{-\alpha}} \right)^{\mult \alpha} = \prod_{n \geq 1} \left( \frac{1-q^n}{1-tq^n} \right)^{\mult n\delta}$, where we let $q:=e^{-\delta}$ throughout. 

Now assume $\kma$ is a twisted affine Kac-Moody algebra, with normalized invariant form $\form[\cdot]{\cdot}$, root lattice $Q$ and null root $\delta$ \cite[Chap. 8]{kac}. Let 
$\kmaf$ denote the underlying finite dimensional simple Lie algebra, and let $\Qf$ denote its root lattice. We define the lattice $M$ as follows \cite{kac}: $M := \Qf$ if $\kma \neq A_{2l}^{(2)}$ and $M := \integers$-span of $2\alpha/\form[\alpha]{\alpha}, \alpha \in \ro(\kmaf)$ (the coroot lattice of $\kmaf$) if $\kma = A_{2l}^{(2)}$.

 For a dominant integral weight $\lambda$ of $\kma$, let $\sch := \sum_\gamma \dim(L(\lambda)_\gamma) \, e^\gamma$ denote the formal character of $L(\lambda)$.  Let $\Lambda_0$ denote the fundamental weight corresponding to the zeroth (additional) node of the Dynkin diagram of $\kma$. Then, $L(\Lambda_0)$ is the basic representation of $\kma$. We have the following classical result \cite{kac}:
\begin{proposition}
Let $\kma$ be a twisted affine algebra. Then
\be
\item 
The formal character of the basic representation of $\kma$ 
is:
$$ e^{-\s \Lambda_0} \, \sch[\s \Lambda_0] = \alolo[1,q] \, \Theta$$
where $\Theta :=\sum_{\alpha \in M} e^\alpha  q^{\form[\alpha]{\alpha}/2}$
is the theta function of the lattice $M$. 
\item The basic string function is given by 
\beq \label{bsf}
\alolo[1,q] = \prod_{n \geq 1} (1-q^n)^{-\mult n\delta}
\eeq
\ee
\end{proposition}

Next, given a formal sum $\xi = \sum_{\alpha \in Q} c_\alpha e^\alpha$, define \cite{igmgod} the {\em constant term} of $\xi$ to be $\ct(\xi) := \sum_{n \in \integers} c_{n\delta} e^{n\delta}$. 
The following simple fact \cite[(5.8)]{svis-kfp} can be used to compute the 
\tsfs $\alm$:
\beq\label{almeq}
\alm = \chkeim \, \ct(e^{-\mu}\,\sch\, \chke)
\eeq
where $\mu$ is a maximal dominant weight of $L(\lambda)$. 
Putting the above facts together, we obtain  the following lemma.
\begin{lemma}\label{ctlem}
Let $\kma$ be a twisted affine algebra. Then,
\be
\item The \tsf of the basic representation of $\kma$ is given by 
$$\alolo = \alolo[1,q] \chkeim \ct(\chke \, \Theta)$$
\item Further, we have $\alolo[1,q] \chkeim = \dprod_{n \geq 1} (1-tq^n)^{-\mult (n\delta)}$
\ee
\end{lemma}

Now, suppose that $\kma$ is a twisted affine algebra of type $X_N^{(r)}$; here  $X_N$  is a simply laced (A-D-E) Dynkin diagram of finite type with a diagram automorphism $\sigma$ of order $r$ ($r=2$ or $3$). Let $\kl$ denote the finite dimensional simple Lie algebra with Dynkin diagram $X_N$ and  let $\sigma$ also denote the corresponding automorphism of $\kl$. For each $k \in \integers$, let $\kl_k \subset \kl$ be the eigenspace of $\sigma$ for the eigenvalue $\exp(2\pi k i/r)$ (so, $\kl_k = \kl_{k+r}$) . Since $\sigma$ acts diagonalizably on $\kl$, we have a $\integers/r\integers$ gradation:
\beq \label{grade}
\kl = \bigoplus_{j \in \integers/r\integers} \kl_j
\eeq
 If $\csa$ is a Cartan subalgebra of $\kl$, let $\csa_j:=\csa \cap \kl_j$ for all $j \in \integers$. We collect together the important facts about the decomposition \eqref{grade}.
\begin{proposition}\label{fdslaut}
{\rm \cite[Chap. 8]{kac}} With notation as above, we have 
\begin{enumerate}
\item  $\kl_0$ is a simple Lie algebra and $\kl_j$ is an irreducible $\kl_0$-module $\forall j$.
\item $\kl_1 \cong \kl_{-1}$ as $\kl_0$-modules.
\item $\csa_0$ is a Cartan subalgebra of $\kl_0$ and its centralizer in $\kl$ is $\csa$.
\item If $\kma$ is not of type $A_{2l}^{(2)}$, then $\kl_0$ and $\kmaf$ are isomorphic. Further, the highest weight of the $\kl_0$-module $\kl_1$ is the dominant short root $\theta_s$ of $\kl_0$.
\item If $\kma$ is of type $A_{2l}^{(2)}$, then $\kl_0$ is of type $B_l$, while $\kmaf$ is of type $C_l$. Further, the highest weight of $\kl_1$ is $2\theta_s$, where $\theta_s$ is the dominant short root of $\kl_0$.
\end{enumerate}
\end{proposition} 

We denote $l:=\text{rank }\kl_0$, $m:=$the number of short simple roots of $\kl_0$ and let $\theta_l$ (resp. $\theta_s$) be the dominant long (resp. short) root of $\kl_0$. 
\begin{proposition}  {\rm(\cite[Chap. 8]{kac} or \cite{carter})}
Let $\kma$ be a twisted affine algebra of type $X_N^{(r)} \neq A_{2l}^{(2)}$.
Consider the action of the cyclic group generated by the automorphism $\sigma$, on the nodes of the Dynkin diagram of $X_N$.
This has the following properties.
\be
\item Each orbit has cardinality 1 or $r$.
\item The number of orbits equals $l$.
\item The number of orbits of cardinality $r$ is equal to $m$.
\item Thus, $m = \frac{N-l}{r-1}$.
\ee
\end{proposition}

Next, we recall that the untwisted affine algebra $\hat{\mathcal{L}}(\kl):=\complex[z,z^{-1}] \otimes \kl \, \oplus \complex K \, \oplus \complex d$, where $K$ is the central element and $d$ is the degree derivation \cite{kac}. We extend $\sigma$ to an automorphism $\tilde{\sigma}$ of $\hat{\mathcal{L}}(\kl)$ by $K \mapsto K$, $d \mapsto d$ and $z^j \otimes x \mapsto  \exp(-2\pi j i/r) \, z^j \otimes \sigma(x)$ for $j \in \integers, x \in \kl$. The fixed point set of $\tilde{\sigma}$ is the affine Lie algebra $\kma = \kma(X_N^{(r)})$.

We have a natural $\integers$-grading $\kma = \oplus_{j \in \integers} \; \kma_j$ with $\kma_0 = \kl_0 + \complex K + \complex d$ and $\kma_j = z^j \otimes \kl_j$ for $j \neq 0$. We observe that for $j \neq 0$, $\kma_j \cong \kl_j$ is an irreducible $\kl_0$-module and that $\kma_j \cong \kma_k$ when $j \equiv k \pmod{r}$, $j, k \neq 0$. Let $\ex_n$ denote the multiset of {\em generalized exponents} of the $\kl_0$-module $\kl_n$ for $n \in \integers$. The main result of this paper is the following.
\begin{theorem} \label{mainthm}
Let $\kma$ be a twisted affine algebra. The \tsf of the basic representation of $\kma$ is given by $$\alolo = \dprod_{n=1}^\infty \, \dprod_{e \in \ex_n} (1- t^{e+1}\, q^n)^{-1}$$
\end{theorem}

The proof will be given in sections \ref{pfa} and \ref{pfb}. But first, we make some remarks.

\begin{remark}
When $\kma$ is an untwisted simply-laced affine, this result was proved in \cite{svis-kfp}. In this case, the $\kl_0$-modules $\kl_j$ are all isomorphic to the adjoint representation of $\kl_0$. Thus $\ex_n = \ex(\kl_0)$, the set of exponents of $\kl_0$ for all $n$, and we recover Theorem \ref{first}.
\end{remark}

\begin{remark}
The cardinality of $\ex_n$ is the dimension of the zero weight space of $\kma_n$. From proposition \ref{fdslaut}, it follows that $|\ex_n| = \dim (z^n \otimes \csa_{n})$. Since $z^n \otimes \csa_{n}$ is the root space of $\kma$ corresponding to the imaginary root $n \delta$, we deduce that $|\ex_n| = \mult (n\delta)$. Thus, this expression is a $t$-deformation of 
the expression for the basic string function (equation \eqref{bsf}).
\end{remark}

\begin{remark} 
From the explicit description of the Chevalley generators of $\kl_0$ in terms of those of $\kl$ \cite[Chap. 8]{kac}, it is clear that a principal nilpotent element of $\kl_0$ is also a principal nilpotent of $\kl$. This observation, together with proposition \ref{fdslaut} implies the following equality of multisets: $$\ex(\kl) = \bigsqcup_{j=1}^r \ex_j$$ where the left hand side is the multiset of exponents of the Lie algebra $\kl$, i.e, the generalized exponents of its adjoint representation. Further, since $\kl_r = \kl_0$, we have $\ex_r = \ex(\kl_0)$. Thus, the sets $\ex(\kl)$ and $\ex(\kl_0
)$  determine the $\ex_n$ for all $n$; this is clear for $r=2$, while for $r=3$ it follows from the further fact that $\ex_1 = \ex_2$.  Table \ref{exptable} lists the $\ex_n$  for all twisted affine algebras.
\end{remark}

\begin{table}
\begin{center}
\setlength{\extrarowheight}{6pt}
 \begin{tabular}{ | c | c | c|}
\hline
$\kma$ & $\ex_0$ & $\ex_{\s 1} = \ex_{\s -1}$\\\hline
$A_{\s 2l}^{(2)}$ ($l \geq 1$) & $1,3,5,\cdots,2l-1$ & $2,4,6,\cdots,2l$ \\\hline
$A_{\s 2l-1}^{(2)}$ ($l \geq 3$) & $1,3,5,\cdots,2l-1$ & $2,4,6,\cdots,2l-2$ \\\hline
$D_{\s l+1}^{(2)}$ ($l \geq 2$) & $1,3,5,\cdots,2l-1$ & $l$ \\\hline
$E_{\s 6}^{(2)}$ & $1, 5, 7, 11$ & $4, 8$ \\ \hline
$D_{\s 4}^{(3)}$ & $1,5$ & $3$ \\\hline
\end{tabular}
\medskip
\caption{\protect $\ex_n$ for the twisted affines $\kma = X_N^{(r)} \; (\ex_{n+r} = \ex_n \text{ for all } n)$.}\label{exptable}
\end{center}
\end{table}

\subsection{}
We derive an interesting corollary of theorem \ref{mainthm}. If $\kma$ is an affine Kac-Moody algebra of rank $l+1$, and $e_i, f_i \, (i=0, \cdots, l)$ are the Chevalley generators, the principal Heisenberg subalgebra $\ph$ of $\kma$ is defined to be $$\ph:=\{x \in \kma: [x, \sum_{i=0}^l \, e_i] \in \complex K\}$$ where $K$ is the central element of $\kma$ \cite{kac}. The principal gradation of $\kma$ induces a gradation $\ph = \oplus_{j \in \integers} \; \ph_j$. If $\kma$ is an untwisted simply-laced or twisted affine algebra, the basic representation $\bas$, as an $\ph$-module, is irreducible. The {\em exponents} of the affine algebra $\kma$ are the elements of the (infinite) multiset $\ex(\kma)$ of nonzero integers in which each $j$ occurs $\dim \ph_j$ times. Let $\ex^+(\kma) := \ex(\kma) \cap \integers_{>0}$ denote the positive exponents of $\kma$. The following lemma relates the multisets $\ex^+(\kma)$ and 
$\ex_n$.

\begin{lemma} \label{affexpo}
Let $\kma$ be a twisted affine algebra or an untwisted simply-laced affine algebra of type $X_N^{(r)}$, with Coxeter number $h$. Then 
$$ \ex^+(\kma) = \{ e + hn: n \geq 0, e \in \ex_n\}$$
\end{lemma}
\smallskip
\noindent
{\em Proof:} Follows easily from \cite[Chap. 14]{kac} and table \ref{exptable}. \qed

\medskip
We deduce the following nice formula for the specialization of the \tsf $\alolo$ at $t \mapsto q, q \mapsto q^h$.
\begin{corollary} \label{spec}
Let $\kma$ be a twisted affine algebra or an untwisted simply-laced affine algebra, with Coxeter number $h$. Let $\kmaf$ be its underlying finite dimensional simple Lie algebra.
Then 
\beq \label{aloeq}
\alolo[q,q^h] = \frac{\dprod_{\bar{e} \in \ex(\kmaf)} (1-q^{\bar{e}+1})}{\dprod_{e \in \ex^{{}^+}\!(\kma)} (1-q^{e+1})}
\eeq
where $\ex(\kmaf)$ is the (finite) multiset of exponents of $\kmaf$. 
\end{corollary}

\smallskip
\noindent
{\em Proof:} Applying the specialization $t \mapsto q, q \mapsto q^h$ to theorem \ref{mainthm}, and using lemma \ref{affexpo}, we obtain equation \eqref{aloeq}, but with $\kl_0$ in place of $\kmaf$. Proposition \ref{fdslaut} implies that $\kl_0$ and $\kmaf$ are either isomorphic or dual. Since dual algebras have the same exponents, the result follows in all cases. \qed

\section{Proof of Theorem \ref{mainthm} for $\kma \neq A_{2l}^{(2)}$} \label{pfa}
Throughout this section, we take $\kma$ to be a twisted affine algebra, $\kma \neq A_{2l}^{(2)}$. Let $\form[\cdot]{\cdot}$ denote the normalized invariant form of $\kma$ \cite{kac}. We then have $\form[\alpha]{\alpha}=2$ for all short real roots of $\kma$. We recall that the {\em height} of a root $\alpha$ (written $\htt \alpha$) is the sum of the coefficients obtained when $\alpha$ is written as a linear combination of simple roots. 
The following result is a special case of Cherednik's {\em difference Macdonald-Mehta} constant term identity \cite{dmmc}.
\begin{proposition} \label{propone}
Let $\kma$ be a twisted affine algebra, $\kma \neq A_{2l}^{(2)}$. Let 
$\kmaf$ be the underlying finite dimensional simple Lie algebra and let $\form[\cdot]{\cdot}$ denote the normalized invariant form of $\kma$. Then we have 
\begin{equation}\label{cmm-id}
\ct(\chke \, \Theta) = \prod_{\alpha \in \rp(\kmaf)} \, \prod_{j=1}^\infty \left( \frac{1 - t^{\htt \alpha} \, q^{\frac{\form[\alpha]{\alpha}}{2} j}}{1 - t^{\htt \alpha + 1} \, q^{\frac{\form[\alpha]{\alpha}}{2} j}} \right)
\end{equation}
\end{proposition}

\pfbegin
The positive real roots of $\kma$ are given by $\rrep(\kma) =\{ \beta + \frac{\form[\beta]{\beta}}{2} j\delta: j \geq 1, \beta \in \ro(\kmaf)\} \cup \rp(\kmaf)$. Thus the Cherednik kernel of $\kma$ becomes  
$$\chke = \prod_{\beta \in \rp(\kmaf)} \, \prod_{j\geq 0} \frac{ (1- e^{-\beta} q^{j\form[\beta]{\beta}/2}) (1- e^{\beta} q^{(j+1)\form[\beta]{\beta}/2}) }{ (1- t e^{-\beta} q^{j\form[\beta]{\beta}/2}) (1- t e^{\beta} q^{(j+1)\form[\beta]{\beta}/2})}$$
Applying \cite[theorem 5.3]{dmmc} with $R$ chosen to be the coroot system of $\kmaf$ yields equation \eqref{cmm-id}. \qed

To simplify notation, we let $\poch[a_1, a_2, \cdots, a_p]{x} := \prod_{i=1}^p \prod_{n=0}^\infty (1 - a_ix^n)$.
Let us now separate the contributions of long and short roots in equation \eqref{cmm-id}. Define 
$$K_s(q) \;\;\;(\text{resp. } K_l(q)) := \prod_{\substack{\alpha \in \rp(\kmaf)\\ \alpha \text{ short} \\ \text{(resp. long)}}} \frac{\poch[t^{\htt \alpha}q]{q}}{\poch[t^{\htt \alpha + 1}q]{q}}$$
Since $\form[\alpha]{\alpha}/2$ is 1 (resp. $r$) if $\alpha$ is short (resp. long), proposition \ref{propone} implies
$$ \ct (\chke \Theta) = K_s(q) K_l(q^r)$$

Now,  for each $k \geq 1$, let $n_k$ (resp. $n_k(s)$) denote the number of positive roots (resp. short positive roots) of $\kmaf$ of height $k$. This gives 
$$ K(q) := K_s(q) K_l(q) = \frac{\poch[tq]{q}^l}{\dprod_{p \geq 1} \poch[t^{p+1}q]{q}^{n_p - n_{p+1}}}$$
where $l$ is the number of simple roots of $\kmaf$. Similarly,
$$ K_s(q) = \frac{\poch[tq]{q}^m}{\dprod_{p \geq 1} \poch[t^{p+1}q]{q}^{n_p(s) - n_{p+1}(s)}}$$
where $m$ is the number of short simple roots of $\kmaf$. We recall the following classical result (see, for example, \cite{ion}) relating the $n_k$ and $n_k(s)$ to generalized exponents of certain representations of $\kmaf$.
\begin{proposition}\label{propion}
With notation as above, $n_p - n_{p+1}$ is the number of times $p$ occurs as an exponent of $\kmaf$ (i.e, as a generalized exponent of the adjoint representation $V(\theta_l)$). Similarly, $n_p(s) - n_{p+1}(s)$ is the number of times $p$ occurs  as a generalized exponent of the representation $V(\theta_s)$ of $\kmaf$.
\end{proposition}

Now, rewriting $K_s(q) K_l(q^r) = K(q^r) \frac{K_s(q)}{K_s(q^r)}$ and using propositions \ref{fdslaut} and \ref{propion}, we get
$$\ct (\chke \Theta) = \poch[tq^r]{q^r}^l \frac{\poch[tq]{q}^m}{\poch[tq^r]{q^r}^m} \dprod_{p \in \ex_0} \frac{1}{\poch[t^{p+1}q^r]{q^r}} \dprod_{p \in \ex_1} \frac{\poch[t^{p+1}q^r]{q^r}}{\poch[t^{p+1}q]{q}}$$
We now observe that $l = \mult j\delta$ for $j \equiv 0 \pmod{r}$ and $m = \frac{N-l}{r-1} = \mult j\delta$ for $j \not\!\equiv 0 \pmod{r}$. Thus, the above equation can be rewritten as :
$$\ct (\chke \Theta) = \dprod_{n \geq 1} (1-tq^n)^{\mult (n\delta)} \prod_{n \geq 1} \prod_{e \in \ex_n} \frac{1}{ 1- t^{e+1}q^n}$$
 Lemma \ref{ctlem} now completes the proof of theorem \ref{mainthm} for all twisted affine algebras $\kma \neq A_{2l}^{(2)}$.

\section{$\kma = A_{2l}^{(2)}$} \label{pfb}
\subsection{}
In this section, we consider the case $\kma = A_{2l}^{(2)}$. The underlying finite dimensional simple Lie algebra is $\kmaf = C_l$. Let $\rof$ be the set of roots of $\kmaf$. Letting $\form[\cdot]{\cdot}$ denote the standard inner product in $\reals^l$ and $\epsilon_i \,(1 \leq i \leq l)$ be the standard orthonormal basis, we can take $\rof = \{\pm \epsilon_i \pm \epsilon_j: 1 \leq i < j \leq l\} \cup \{\pm 2\epsilon_i: 1 \leq i \leq l\}$. We observe that the coroot lattice $M$ of $\kmaf$ is just $M = \oplus_{i=1}^l \integers \epsilon_i$. The set of real roots of $\kma$ is given by $\rre = S_1 \cup S_2 \cup S_4$ where $S_1 = \{\frac{1}{2} (\alpha + (2n-1)\delta): \alpha \in \rof_l\}$, $S_2 = \{\alpha + n \delta: \alpha \in \rof_s\}$ and $S_4 = \{\alpha + 2n\delta: \alpha \in \rof_l\}$, where $\delta$ is the null root of $\kma$ and $\rof_l$ (resp. $\rof_s$) denotes the set of long (resp. short) roots in $\rof$. The elements of $S_n$ have norm $n$ ($n=1, 2, 4$), and each $S_n$ is invariant under the Weyl group $W$ of $\kma$. Let $\chke$ denote the Cherednik kernel of $A_{2l}^{(2)}$, given by equation \eqref{which}.

Now enlarge $\rre$ by defining:
$$\cc := \bigcup_{i=1}^5 \cc_i \text{ where } \cc_1 := (1/2)S_4, \,\cc_2 := S_4, \,\cc_3 := S_1, \,\cc_4:= 2S_1, \cc_5:= S_2$$ The set $\cc$ is the non-reduced irreducible affine root system of type $\ccheckc$ in the classification of Macdonald \cite{aha-bk}. Observe that $\cc$ is $W$-invariant, with each $\cc_i$ being a $W$-orbit. Following the notation of Macdonald \cite{aha-bk}, define $R_1^+:=\{ \epsilon_1, \cdots, \epsilon_l\}$ and $R_2^+:=\{\epsilon_i \pm \epsilon_j: 1 \leq i < j \leq l\}$. We now let $k_i$ ($1 \leq i \leq  5$) be arbitrary parameters, and let $(u_1, u_2, u_3, u_4) = (q^{k_1}, -q^{k_2},q^{k_3+\frac{1}{2}}, - q^{k_4+\frac{1}{2}})$ and $(u'_1, u'_2, u^\prime_3, u^\prime_4) = (qu_1, qu_2, u_3, u_4)$. The Cherednik kernel $\spchke$ (with parameters $k_i$) for the non-reduced affine root system $\cc$ then becomes $\spchke : = \spchke^{(1)} \spchke^{(2)}$ where
\begin{align*}
\spchke^{(1)} &= \dprod_{\alpha \in R_1^+} \frac{\poch[e^{-2\alpha}, qe^{2\alpha}]{q}}{\prod_{i=1}^4 \poch[u_ie^{-\alpha}]{q} \poch[u^\prime_ie^{\alpha}]{q}}\\
\spchke^{(2)} &= \dprod_{\alpha \in R_2^+} \frac{\poch[e^{-\alpha}, qe^{\alpha}]{q}}{\poch[q^{k_5}e^{-\alpha}, q^{k_5+1}e^{\alpha}]{q}}
\end{align*}
\cite[(5.1.14)]{aha-bk}. The following lemma relates the kernels $\spchke$ and $\chke$.

\begin{lemma}
Define $t := q^{k_5}$, and let the parameters $k_i$ satisfy the relations $k_3 = k_5 = 2k_1 = 2k_2$. We then have:
\be
\item If $k_4=0$, then $\spchke = \chke$.
\item If $k_4 \to \infty$ (i.e $q^{k_4} \to 0$), then $\spchke \to \dfrac{\chke \, \Theta_M}{\poch[q]{q}^l}$ where $M :=\oplus_{i=1}^l \integers \epsilon_i$ and $\Theta_M :=\sum_{\alpha \in M} e^\alpha  q^{\form[\alpha]{\alpha}/2}$ is its theta function.
\ee
\end{lemma}
\pfbegin The first statement is easy; in fact one can recover the Cherednik kernels of all twisted affines (and all non-reduced affine root sytems) by appropriate specialization of $\spchke$ \cite[(5.1.7)]{aha-bk}. To prove (2), we observe that for the given choice of parameters, one gets $$\spchke =  \chke \, \prod_{i=1}^l \poch[-q^{1/2} e^{\epsilon_i},-q^{1/2} e^{-\epsilon_i}]{q}$$ 
Now, by the Jacobi triple product identity, we have 
$$\poch[-q^{1/2} e^{\epsilon_i},-q^{1/2} e^{-\epsilon_i}]{q} = \poch[q]{q}^{-1} \, \sum_{n \in \integers} q^{n^2/2} e^{n \epsilon_i} = \poch[q]{q}^{-1} \,\Theta_{\integers \epsilon_i}$$
Since the theta function $\Theta_M$ of the rectangular lattice $M =\oplus_{i=1}^l \integers \epsilon_i$ is just the product $\prod_{i=1}^l \Theta_{\integers \epsilon_i}$, the result follows. \qed

To obtain the explicit form of the Cherednik-Macdonald-Mehta constant term identity for $A_{2l}^{(2)}$ (i.e an expression for $\ct(\chke \, \Theta_M)$) it thus only remains to specialize the well-known formula for $\ct (\spchke)$ \cite[(5.8.20)]{aha-bk} at $k_3 = k_5 = 2k_1 = 2k_2$ and $k_4 \to \infty$ (and letting $q^{k_5}=:t$). We give the result of this (long, but straightforward) calculation below.

\begin{proposition}\label{cmm-att}
For $\kma = A_{2l}^{(2)}$, we have:
$$\ct(\chke \, \Theta_M) = \dfrac{\poch[tq]{q}^l}{\poch[t^2q^2, t^4q^2, \cdots, t^{2l}q^2]{q^2}\poch[t^3q, t^5q, \cdots, t^{2l+1}q]{q^2}}$$
\end{proposition} 
Now, for $A_{2l}^{(2)}$, recall that $\mult (j \delta) = l \, \forall j \geq 1$. This fact, together with proposition \ref{cmm-att}, lemma \ref{ctlem} and Table \ref{exptable} proves theorem \ref{mainthm} for $\kma = A_{2l}^{(2)}$.

\subsection{Two variable generalization} One can prove a slightly more general, two-variable version of theorem \ref{mainthm} for $\kma = A_{2l}^{(2)}$. To state this, let $s, t$ be indeterminates, and define the two-variable Kostant partition function $\stko[\beta]$ to be the coefficient of $e^{\beta}$  in the product 
$\dprod_{\alpha \in \rp(\kma)} (1-u_\alpha \, e^{\alpha})^{-\mult \alpha}$ where $u_\alpha:=s$ if $\alpha$ is a real root of norm 1 (=shortest root length) and $u_\alpha:=t$ for all other roots (i.e imaginary roots, and real roots of norms 2 and 4). 
For a dominant integral weight $\lambda$ of $A_{2l}^{(2)}$, and a maximal dominant weight $\mu$ of $L(\lambda)$, 
define the two variable Kostka-Foulkes polynomial 
$$\kostkast[\lambda\mu] := \sum_{w \in W} \epsilon(w)\, \stko[w(\lambda + \rho) - (\mu + \rho)]$$ 
and let the corresponding $(s,t)$-string function for the basic representation be 
$\alolost := \displaystyle\sum_{\s p \geq 0} \kostkast[\s \Lambda_0, \s \Lambda_0 - p \delta] \, q^p$. The following is the two variable version of theorem \ref{mainthm}.
\begin{proposition}
For $\kma =  A_{2l}^{(2)}$, $$\alolost = \prod_{\substack{j=1\\j \text{ even}}}^{2l} \poch[t^jq^2]{q^2}^{-1} \, \prod_{\substack{j=1\\j \text{ odd}}}^{2l} \poch[s^2t^jq]{q^2}^{-1}$$
\end{proposition}

\pfbegin
The proof is along the exact same lines as that of proposition \ref{cmm-att}, but now with parameters chosen differently. We choose $k_5 = 2k_1 = 2k_2$, $k_4 \to \infty$, but leave $k_3$ as a free parameter. We then take $t:=q^{k_5}$ and $s:=q^{k_3}$. The remaining details are easily checked. \qed
 
\begin{corollary}\label{posi}
$\kostkast[\s \Lambda_0, \s \Lambda_0 - p \delta] \in \integers_{\geq 0}[s,t]$ for all $p \geq 0$.
\end{corollary}

Finally, we remark that it would be of interest to find a more natural explanation for the positivity result of the above corollary (or more generally, for $\kostkast[\lambda\mu]$) in terms of a Brylinski-Kostant type filtration, as is known for the usual (one variable) affine Kostka-Foulkes polynomials \cite{slofstra}. We also note that the two variable Kostka-Foulkes polynomials can be defined for all twisted affines (in fact, for any affine root system with more than one root length) and in more than one way (corresponding to different choices of the $u_\alpha$ in the definition). But it appears, from preliminary calculations, that only $A_{2l}^{(2)}$ (with the given choice of $u_\alpha$) exhibits the positivity property of corollary \ref{posi}.

%% file: twisted-main.bbl
\begin{thebibliography}{10}

\bibitem{carter}
R.~W. Carter.
\newblock {\em Lie algebras of finite and affine type}, volume~96 of {\em
  Cambridge Studies in Advanced Mathematics}.
\newblock Cambridge University Press, Cambridge, 2005.

\bibitem{dmmc}
I.~Cherednik.
\newblock Difference {M}acdonald-{M}ehta conjecture.
\newblock {\em Internat. Math. Res. Notices}, 10:449--467, 1997.

\bibitem{frenkel-kac}
I.~B. Frenkel and V.~G. Kac.
\newblock Basic representations of affine {L}ie algebras and dual resonance
  models.
\newblock {\em Invent. Math.}, 62(1):23--66, 1980/81.

\bibitem{rkg}
R.~K. Gupta.
\newblock Characters and the {$q$}-analog of weight multiplicity.
\newblock {\em J. London Math. Soc. (2)}, 36(1):68--76, 1987.

\bibitem{ion}
B.~Ion.
\newblock The {C}herednik kernel and generalized exponents.
\newblock {\em Int. Math. Res. Not.}, (36):1869--1895, 2004.

\bibitem{kac}
V.~G. Kac.
\newblock {\em Infinite dimensional {L}ie algebras}.
\newblock Cambridge University Press, third edition, 1990.

\bibitem{kac-kazhdan-lepow-wilson}
V.~G. Kac, D.~A. Kazhdan, J.~Lepowsky, and R.~L. Wilson.
\newblock Realization of the basic representations of the {E}uclidean {L}ie
  algebras.
\newblock {\em Adv. in Math.}, 42(1):83--112, 1981.

\bibitem{kacpeterson}
V.~G. Kac and D.~H. Peterson.
\newblock Infinite-dimensional {L}ie algebras, theta functions and modular
  forms.
\newblock {\em Adv. in Math.}, 53(2):125--264, 1984.

\bibitem{kac-peterson-112constr}
V.~G. Kac and D.~H. Peterson.
\newblock {$112$} constructions of the basic representation of the loop group
  of {$E_8$}.
\newblock In {\em Symposium on anomalies, geometry, topology ({C}hicago,
  {I}ll., 1985)}, pages 276--298. World Sci. Publishing, Singapore, 1985.

\bibitem{kostant3d}
B.~Kostant.
\newblock The principal three-dimensional subgroup and the {B}etti numbers of a
  complex simple {L}ie group.
\newblock {\em Amer. J. Math.}, 81:973--1032, 1959.

\bibitem{aha-bk}
I.~G. Macdonald.
\newblock {\em Affine {H}ecke algebras and orthogonal polynomials}, volume 157
  of {\em Cambridge Tracts in Mathematics}.
\newblock Cambridge University Press, Cambridge, 2003.

\bibitem{igmgod}
I.~G. Macdonald.
\newblock A formal identity for affine root systems.
\newblock In {\em Lie groups and symmetric spaces}, volume 210 of {\em Amer.
  Math. Soc. Transl. Ser. 2}, pages 195--211. Amer. Math. Soc., Providence, RI,
  2003.

\bibitem{announcement}
Sachin~S. Sharma and Sankaran Viswanath.
\newblock The $t$-analog of the level one string function for twisted affine
  {K}ac-{M}oody algebras.
\newblock {\em Comptes Rendus Mathematique}, 2012.
\newblock to appear.

\bibitem{slofstra}
W.~Slofstra.
\newblock A {B}rylinski filtration for affine {K}ac–-{M}oody algebras.
\newblock {\em Advances in Mathematics}, 229(2):968 -- 983, 2012.

\bibitem{svis-kfp}
S.~Viswanath.
\newblock Kostka-{F}oulkes polynomials for symmetrizable {K}ac-{M}oody
  algebras.
\newblock {\em S\'em. Lothar. Combin.}, 58:Art. B58f, 2008.

\end{thebibliography}
